\documentclass[reqno]{amsart}
\usepackage{latexsym,amscd,amssymb}
\newtheorem{thm}{Theorem}
\newtheorem{lem}{Lemma}

\theoremstyle{definition}

\theoremstyle{remark}

\newtheorem{remark}{Remark}

\newcommand{\I}{\mathcal{I}}

\newcommand{\K}{\mathcal{K}}
\newcommand{\M}{\mathcal{M}}

\newcommand{\Q}{\mathbb{Q}}
\newcommand{\T}{\mathcal{T}}

\newcommand{\siegel}{\mathfrak{S}}
\newcommand{\Z}{\mathbb{Z}}

\newcommand{\re}{{\operatorname{Re}}}

\begin{document}

\title[Homology of Torelli groups]{Homological infiniteness of \\
Torelli groups}
\author{Toshiyuki Akita}
\address{Department of Applied Mathematics, Fukuoka University,
Fukuoka 814-80, Japan}
\email{akita@sm.fukuoka-u.ac.jp}
\thanks{The author is supported by Grand-in-Aid for
Encouragement of Young Scientists (No. 09740072),
the Ministry of Education, Science,
Sports and Culture.}


\maketitle
\section{Introduction}
Let $\Sigma_g$ be a closed orientable surface of genus $g\geq 2$.
Let $\M_{g,r}^n$ be the mapping class group of $\Sigma_g$ relative to
$n$ distinguished points and $r$ fixed embedded disks.
The action of  $\M_{g,r}^n$ on the homology of $\Sigma_g$
induces a surjective homomorphism
\[
\M_{g,r}^n\rightarrow Sp(2g,\Z),
\]
where $Sp(2g,\Z)$ is the Siegel modular group.
The Torelli group $\I_{g,r}^n$ is defined to be its kernel so that we
have an extension
\[
1\rightarrow\I_{g,r}^n\rightarrow\M_{g,r}^n\rightarrow
Sp(2g,\Z)\rightarrow 1.
\]
We omit the decorations $n$ and $r$ when they are zero.

In a series of papers
\cite{johnson0,johnson1,johnson2,johnson3},
D. Johnson obtained several fundamental results
concerning the structure of $\I_g$ and $\I_{g,1}$ (see also
\cite{johnson-survey,hain-survey}).
In particular, he proved that
$\I_g$ and $\I_{g,1}$ are finitely generated for all $g\geq 3$.
On the contrary, A. Miller and D. McCullough \cite{miller} showed that
$\I_2$
(and hence $\I_{2,r}^n$ for all $n,r\geq0$)
is not finitely generated.
G. Mess \cite{mess} showed that $\I_2$ is a free group on
infinitely many generators.
Johnson and J. Millson showed that $H_3(\I_3,\Z)$ contains a free
abelian group of infinite rank (cf. \cite{mess}).
It is not known whether $\I_{g,r}^n$ is finitely presented for $g\geq
3$.
In this paper, we will prove:
\begin{thm}\label{thm-torelli}
For all $n\geq 0$, the rational homology $H_*(\I_g^n,\Q)$ of the
Torelli group $\I_g^n$ is infinite dimensional over $\Q$ if $g$ is
sufficiently large compared with $n$.
In particular, $H_*(\I_g,\Q)$ and $H_*(\I_{g}^1,\Q)$ are infinite
dimensional for $g\geq 7$.
\end{thm}\noindent
This theorem yields the negative answer to the question posed
by Johnson \cite{johnson-survey} which asks whether the Torelli
space $\mathbf{T}_{g}^1$ (see \S 2 for the definition) is homotopy
equivalent to a finite complex, provided $g\geq 7$.

For $n+r\leq 1$, let $\K_{g,r}^n$ be the subgroup of $\M_{g,r}^n$
generated by all the Dehn twists along separating simple closed
curves.
The groups $\K_{g,1}$ and $\K_g$
are related to the Casson invariants of homologly 3-spheres
through the work of S. Morita \cite{morita1989,morita1991}.
For $g=2$, $\K_2$ is equal to $\I_2$ so that it is a free 
group on infinitely many generators.
In contrast, the group $\K_g$ is not free for $g\geq 3$ and
almost nothing is known about the structure of this group, however, we
can prove:
\begin{thm}\label{thm-kg}
For all $g\geq 3$,
$H_*(\K_g,\Q)$ and $H_*(\K_{g}^1,\Q)$ are infinite dimensional over
$\Q$.
\end{thm}\noindent

\begin{remark}
There is no general agreement on the definition of $\I_{g,r}^n$ when
$r+n>1$. We employ the one which was used in \cite{hain-survey}.
It differs from that which was given in \cite{johnson2}.
\end{remark}

\section{Preliminaries}
In this section, we recall relevant definitions and facts concerning
of Torelli groups which will be used later.
The reader should refer to \cite{hain-survey,harer} for further
detail.

Let $\T_g^n$ be the Teichm\"uller space with $n$ marked points.
$\M_g^n$ acts on $\T_g^n$ properly discontinuously and the quotient
space $\mathbf{M}_g^n=\T_g^n/\M_g^n$ is, by definition,  the moduli
space of curves of genus $g$ with $n$ marked points.
The Torelli group $\I_g^n$ is torsion-free and hence it acts on
$\T_g^n$ freely so that the quotient space
$\mathbf{T}_g^n=\T_g^n/\I_g^n$ is a complex manifold.
Moreover, since $\T_g^n$ is contractible, 
$\mathbf{T}_g^n$ is the classifying space of $\I_g^n$ so that
there is a canonical isomorphism
\[
H_*(\I_{g}^n,\Z)\cong H_*(\mathbf{T}_{g}^n,\Z).
\]
$\mathbf{T}_g^n$ is called the {\em Torelli space} and is important in 
algebraic geometry (cf. \cite{hain-survey}).
The action of $\M_g^n$ on $\T_g^n$ induces the properly
discontinuous action of $Sp(2g,\Z)=\M_g^n/\I_g^n$ on $\mathbf{T}_g^n$ so
that the quotient space $\mathbf{T}_g^n/Sp(2g,\Z)$ coincides with
$\mathbf{M}_g^n$.

Let $\siegel_g$ be the Siegel upper half space of degree $g$.
The group $Sp(2g,\Z)$ acts on $\siegel_g$ properly discontinuously and 
the quotient space $Sp(2g,\Z)\backslash\siegel_g$ is identified with
the moduli space of principally polarized abelian varieties of
dimension $g$.
For an integer $L\geq 3$,
let $\Gamma(L)$ be the principal congruence subgroup of $Sp(2g,\Z)$ of
level $L$ defined to be the kernel of the canonical homomorphism
\[
Sp(2g,\Z)\rightarrow Sp(2g,\Z/L\Z).
\]
$\Gamma(L)$ is a torsion-free subgroup of finite index of $Sp(2g,\Z)$
and hence acts on $\siegel_g$ freely.
The quotient space $\Gamma(L)\backslash\siegel_g$ is identified with
the moduli space of principally polarized abelian varieties of
dimension $g$ with level $L$ structure.
It is known that $\Gamma(L)\backslash\siegel_g$ is homotopy
equivalent to a finite complex (see
\cite{borel-serre,serre-arithmetic} for instance).

Let $\M_g^n(L)\subset\M_g^n$ be the full inverse image of $\Gamma(L)$
under the homomorphism $\M_g^n\rightarrow Sp(2g,\Z)$ so that it fits
into the extension
\begin{equation}\label{mg-level-l}
1\rightarrow \I_g^n\rightarrow \M_g^n(L)\rightarrow
\Gamma(L)\rightarrow 1.
\end{equation}
$\M_g^n(L)$ is a torsion-free subgroup of finite index of $\M_g^n$ and
hence acts on $\T_g^n$ freely.
The quotient space $\mathbf{M}_g^n(L)=\T_g^n/\M_g^n(L)$ is, by
definition,  the moduli space of curves of genus $g$ with $n$ marked
points and level $L$ structure.
The action of $\M_g^n(L)$ on $\T_g^n$ induces the free action
of $\Gamma(L)=\M_g^n(L)/\I_g^n$ on $\mathbf{T}_g^n$ so that the
quotient space $\mathbf{T}_g^n/\Gamma(L)$ coincides with
$\mathbf{M}_g^n(L)$.
According to the work of W. Harvey \cite{harvey}, $\mathbf{M}_g^n(L)$
is homotopy equivalent to a finite complex.

\section{Proof of Theorem \ref{thm-torelli}}
Fix an integer $L\geq 3$.
Since $\siegel_g$ is contractible, the projection
$\siegel_g\rightarrow\Gamma(L)\backslash\siegel_g$ is the
universal principal $\Gamma(L)$-bundle.
The associated bundle
\begin{equation}\label{borel-const}
\mathbf{T}_g^n\rightarrow\siegel_g\times_{\Gamma(L)}\mathbf{T}_g^n
\rightarrow \Gamma(L)\backslash\siegel_g.
\end{equation}
is nothing but the Borel construction of the $\Gamma(L)$-space
$\mathbf{T}_g^n$.
Since $\Gamma(L)$ acts freely on $\mathbf{T}_g^n$, the total space
$\siegel_g\times_{\Gamma(L)}\mathbf{T}_g^n$ is homotopy equivalent to
$\mathbf{T}_g^n/\Gamma(L)=\mathbf{M}_g^n(L)$.
Note that the associated bundle (\ref{borel-const}) is identified, up
to homotopy, with the fibration $B\I_g^n\rightarrow
B\M_g^n(L)\rightarrow B\Gamma(L)$ of classifying spaces induced from
the exact sequence (\ref{mg-level-l}).

Now suppose that $H_*(\I_g^n,\Q)\cong H_*(\mathbf{T}_g^n,\Q)$ is
finite dimensional.
As $\Gamma(L)\backslash\siegel_g$ is homotopy equivalent to a finite
complex, we may apply the following lemma to the associated bundle
(\ref{borel-const}):
\begin{lem}\label{lem-q-euler-char}
Let $F\rightarrow E\rightarrow B$ be a fibration such that $B$ is a
finite complex and $\dim_{\Q}H_*(F,\Q)<\infty$.
Then $\dim_{\Q}H_*(E,\Q)<\infty$ and
\[
\chi_{\Q}(E)=\chi_{\Q}(F)\cdot\chi(B),
\]
where $\chi_{\Q}$ is defined by
$\chi_{\Q}(-)=\sum_i(-1)^i\dim_{\Q}H_i(-,\Q)$.
\end{lem}\noindent
This lemma is a direct consequence of the Serre spectral sequence
applied to the fibration $F\rightarrow E\rightarrow B$.
As a result, one has
\[
\chi(\mathbf{M}_g^n(L))
=\chi_{\Q}(\mathbf{T}_g^n)\cdot\chi(\Gamma(L)\backslash\siegel_g).
\]
Since both of the projections
$\Gamma(L)\backslash\siegel_g\rightarrow Sp(2g,\Z)\backslash\siegel_g$
and $\mathbf{M}_g^n(L)\rightarrow\mathbf{M}_g^n$ are
$|Sp(2g,\Z/L\Z)|$-fold branched coverings, one has
\begin{align*}\label{prod-formula}
 \chi(\Gamma(L)\backslash\siegel_g) &=|Sp(2g,\Z/L\Z)|\cdot
e(Sp(2g,\Z)\backslash\siegel_g) \\
 \chi(\mathbf{M}_g^n(L)) &=|Sp(2g,\Z/L\Z)|\cdot
e(\mathbf{M}_g^n),
\end{align*}
and hence
\begin{equation}\label{prod-formula}
e(\mathbf{M}_g^n)=\chi_{\Q}(\mathbf{T}_g^n)\cdot
e(Sp(2g,\Z)\backslash\siegel_g),
\end{equation}
where $e$ denotes the orbifold Euler characteristics.
According to G. Harder \cite{harder}, one has
\[
e(Sp(2g,\Z)\backslash\siegel_g)=\prod_{k=1}^g\zeta(1-2k),
\]
while according to J. Harer and D. Zagier \cite{harer-zagier}, one has 
\[
e(\mathbf{M}_g^n)=
\begin{cases}{\displaystyle\frac{1}{2-2g}\zeta(1-2g)}&\mbox{ if }n=0\\
{\displaystyle (-1)^{n-1}\frac{(2g+n-3)!}{(2g-2)!}\zeta(1-2g)}
&\mbox{ if }n>0,
\end{cases}
\]
where $\zeta$ is the Riemman $\zeta$-function.
See also \cite{penner,kont}.
Hence the equality (\ref{prod-formula}) leads to
\[
\chi_{\Q}(\mathbf{T}_g^n)=
\begin{cases}
{\displaystyle\frac{1}{2-2g}\prod_{k=1}^{g-1}
\frac{1}{\zeta(1-2k)}}&\mbox{ if }n=0\\
{\displaystyle (-1)^{n-1}\frac{(2g+n-3)!}{(2g-2)!}
\prod_{k=1}^{g-1}\frac{1}{\zeta(1-2k)}}&\mbox{ if }n>0
\end{cases}
\]
By definition, $\chi_{\Q}(\mathbf{T}_g^n)$ must be an integer and
the proof of Theorem \ref{thm-torelli} is then completed by
virtue of the following lemma which will be proven in the next
section.
\begin{lem}\label{lem-zeta}
For positive integers $m,n$, set
\[
e(m,n)= 
\frac{(2m+n-1)!}{(2m)!}
\prod_{k=1}^{m}\frac{1}{|\zeta(1-2k)|}.
\]
Then, for each $n\geq 1$, $e(m,n)$ is not an integer for sufficiently large
$m$ compared with $n$.
In particular, $e(m,1)$ is not an integer for all $m\geq 6$.
\end{lem}\noindent

\section{The proof of Lemma \ref{lem-zeta}}
Recall that the Riemann $\zeta$-function is defined for
$\re\ s>1$ by
\[
\zeta(s)=\sum_{n=1}^{\infty}\frac{1}{n^s},
\]
On the other hand, the equality
\[
\zeta(1-2k)=(-1)^k\frac{2\cdot (2k-1)!}{(2\pi)^{2k}}\zeta(2k).
\]
holds for any integer $k\geq 1$.
It follows that
\[ 
|\zeta(1-2k)| 
=\frac{2\cdot (2k-1)!}{(2\pi)^{2k}}
\sum_{n=1}^{\infty}\frac{1}{n^{2k}} 
>\frac{2\cdot (2k-1)!}{(2\pi)^{2k}}
\] 
for any integer $k\geq 1$, and hence
\[
e(m,n)<
\frac{(2m+n-1)!}{(2m)!}\prod_{k=1}^{m}
\frac{(2\pi)^{2k}}{2\cdot (2k-1)!}
\]
for any integer $m\geq 1$.
We claim that the right hand side of the inequality converges
to $0$ as $m\rightarrow\infty$.
Indeed, regarding the right hand side as a numerical sequence with
respect to $m$, the ratio of the $(m+1)$-th term to the $m$-th term is 
given by
\[
\frac{(2m+n+1)(2m+n)}{(2m+2)(2m+1)}\cdot
\frac{(2\pi)^{2m+2}}{2\cdot (2m+1)!}.
\]
This converges to $0$ as $m\rightarrow\infty$, hence verifying the
claim.
We conclude that $e(m,n)<1$ and hence $e(m,n)$ is not an integer for
$m$ sufficiently large compared with $n$.
The first assertion is proved.

To prove the second assertion, observe that
${(2\pi)^{2k}}/({2\cdot (2k-1)!})<1$ for $k\geq 9$.
It follows that $e(m,1)$ is strictly decreasing with respect to $m$
for $m \geq 9$.
On the other hand, with the help of a computer, one has
\[
\prod_{k=1}^{14}\zeta(1-2k)=-297203.11\cdots .
\]
We see that, for all $m\geq 14$, $e(m,1)<1$ and hence $e(m,1)$ is
not an integer.
It remains to be proven that $e(m,1)$ is not an integer for $6\leq
m\leq 13$.
However, this can be verified by direct calculations and we omit
the detail.
\begin{remark}
Actually, $e(m,n)$ is not an integer for $n<678$ and $m\geq 6$ and
hence $H_*(\I_g^n,\Q)$ is infinite dimensional for $n<678$ and $g\geq
7$.
We describe briefly how this can be proven.
Recall that $\zeta(1-2k)$ for positive integer $k$ is given by
\[
\zeta(1-2k)=-\frac{B_{2k}}{2k}\in\Q,
\]
where $B_{2k}$ is the $2k$-th Bernoulli number defined as the
coefficient of $z^{2k}/(2k)!$ in the power series expansion of
$z/(e^z-1)$.
On the other hand,  von Staudt's theorem asserts
that the denominator of $B_{2k}$ is not divisible by a prime $p$ if
$2k<p-1$.
By applying von Staudt's theorem to the primes 691 and 3617,
the numerators of $B_{12}$ and $B_{16}$ respectively,
we see that $e(m,n)$ is not an integer for $n<678$ and $6\leq
m\leq1470$.
Now the assertion follows from the inequality
\[
e(m,n)<\frac{(2m+677)!}{(2m)!}\cdot\prod_{k=1}^m
\frac{(2\pi)^{2m+2}}{2\cdot (2m+1)!}<1
\]
which holds for $n<678$ and $m\geq 37$.
%
\end{remark}

\section{Proof of Theorem \ref{thm-kg}}
To prove Theorem \ref{thm-kg}, we first
recall some of Johnson's results concerning the Torelli groups.
Suppose $g\geq 3$ and $[\Sigma_g]\in\wedge^2 H_1(\Sigma_g,\Z)$
corresponds to the fundamental class of $\Sigma_g$.
Under these conditions, Johnson constructed in \cite{johnson0}
natural $Sp(2g,\Z)$-equivariant surjective homomorphisms
\[
\tau_{2,1}:\I_{g,1}\rightarrow \wedge^3 H_1(\Sigma_g,\Z)
\]
and
\[
\tau_{2}:\I_g\rightarrow \wedge^3 H_1(\Sigma_g,\Z)/([\Sigma_g]\wedge
H_1(\Sigma_g,\Z))
\]
and proved in \cite{johnson2} that $\ker\tau_{2,1}=\K_{g,1}$ and
$\ker\tau_{2}=\K_g$.
The homomorphisms $\tau_{2,1}$ and $\tau_2$ are called {\em Johnson
homomorphisms}.
For simplicity, we abbreviate
$\wedge^3 H_1(\Sigma_g,\Z)/([\Sigma_g]\wedge H_1(\Sigma_g,\Z))$
by $\wedge^3 H/H$ and
$\wedge^3 H_1(\Sigma_g,\Z)$ by
$\wedge^3 H$.
As a consequence, $\K_g$ fits into the extension
\[
1\rightarrow\K_g\rightarrow\I_g\stackrel{\tau_2}{\rightarrow}
\wedge^3 H/H\rightarrow 1
\]
Take the classifying space of each group in the extension to yield a
fibration
\begin{equation}\label{fibre-kg}
B\K_{g}\rightarrow B\I_{g}\rightarrow B(\wedge^3 H/H).
\end{equation}
Observe that $B\I_{g}$ is homotopy equivalent to $\mathbf{T}_{g}$ and
$B(\wedge^3 H/H)$ is homotopy equivalent to the
$\binom{2g}{3}-2g$-dimensional torus since $\wedge^3
H/H$ is a free abelian group of rank $\binom{2g}{3}-2g$.

Now suppose $H_*(\K_g,\Q)\cong H_*(B\K_g,\Q)$ is finite dimensional.
It follows from Lemma \ref{lem-q-euler-char} that
$\dim_{\Q}H_*(\mathbf{T}_g,\Q)<\infty$.
If $\dim_{\Q}H_*(\mathbf{T}_g,\Q)<\infty$ (and hence $g\leq 6$), then,
as in the proof of Theorem \ref{thm-torelli},
$\chi_{\Q}(\mathbf{T}_g)$ is defined and satisfies
\[
\chi_{\Q}(\mathbf{T}_g)=\frac{1}{2-2g}\prod_{k=1}^{g-1}
\frac{1}{\zeta(1-2k)}\not=0.
\]
On the other hand,
by applying Lemma \ref{lem-q-euler-char} to the fibration
(\ref{fibre-kg}), one has
$\chi_{\Q}(\mathbf{T}_g)=\chi_{\Q}(B\K_g)\cdot \chi(B(\wedge^3
H/H))=0$ since $\chi(B(\wedge^3 H/H))=0$.
A contradiction.

To prove Theorem \ref{thm-kg} for the group $\K_g^1$, we will identify
$\K_g^1$ with the kernel of a variant of the Johnson homomorphism.
Recall that Torelli groups $\I_g^1$ and $\I_{g,1}$ fit into the
central extension
\[ 
1\rightarrow\Z\rightarrow\I_{g,1}\rightarrow\I_{g}^1\rightarrow 1,
\]
where the center $\Z$ is generated by the Dehn twist $\xi$ along a
simple closed curve parallel to the boundary of a fixed embedded disk
$D\subset\Sigma_g$.

Now the Dehn twist $\xi$ is contained in
$\K_{g,1}$ and hence $\tau_{2,1}$ induces a homomorphism
$\tau_2^1:\I_g^1\rightarrow \wedge^3 H$.
We claim that $\ker\tau_2^1=\K_{g}^1$.
Indeed, $\ker\tau_2^1$ coincides with the image of $\K_{g,1}$ under
the homomorphism $\I_{g,1}\rightarrow\I_g^1$.
But the image of $\K_{g,1}$ is nothing but $\K_g^1$ since any Dehn
twist along separating simple closed curve is isotopic to that
which fixes the embedded disk pointwise.
In summary, we have the extension
\[
1\rightarrow\K_{g}^1\rightarrow\I_{g}^1
\stackrel{\tau_2^1}{\rightarrow}\wedge^3 H\rightarrow 1.
\]
Take the classifying space for each group in the extension to yield a
fibration
\[ 
B\K_{g}^1\rightarrow B\I_{g}^1\rightarrow B(\wedge^3 H).
\]
Now $B\I_{g}^1$ is homotopy equivalent to $\mathbf{T}_{g}^1$ and
$B(\wedge^3 H)$ is homotopy equivalent to the
$\binom{2g}{3}$-dimensional torus since
$\wedge^3 H$ is a free abelian group of rank $\binom{2g}{3}$.
The rest of the proof is similar to that of $\K_g$.
\vspace{3mm}\\
{\em Acknowledgement.}
The author thanks to Professor Shigeyuki Morita for calling author's
attention to the groups $\K_{g,r}^n$.

\ifx\undefined\bysame
\newcommand{\bysame}{\leavevmode\hbox to3em{\hrulefill}\,}
\fi

\end{document}